\documentclass[12pt]{amsart}

\makeatletter
\@namedef{subjclassname@2020}{\textup{2020} Mathematics Subject Classification}
\makeatother

\newcommand{\sech}{\tn{sech}}
\newcommand{\abs}[1]{\ensuremath{\left|{#1}\right|}}

\newcommand{\mg}{\infty}

\newcommand{\pz}{\pi z}
\newcommand{\piz}{\pi z}

\newcommand{\tom}{\to_{M}}

\newcommand{\dsum}{\di\sum}

\newcommand{\hgt}[1]{\rule[0pt]{0pt}{#1}}

\newcommand{\dep}[1]{{\vrule width 0pt height 0pt depth #1}}

\DeclareFontFamily{U}{mathx}{\hyphenchar\font45}
\DeclareFontShape{U}{mathx}{m}{n}{
      <5> <6> <7> <8> <9> <10>
      <10.95> <12> <14.4> <17.28> <20.74> <24.88>
      mathx10
      }{}
\DeclareSymbolFont{mathx}{U}{mathx}{m}{n}
\DeclareFontSubstitution{U}{mathx}{m}{n}
\DeclareMathAccent{\widecheck}{0}{mathx}{"71}

\renewcommand{\kill}[1]{}
\newcommand{\dummy}[1]{\mbox{}}

\makeatletter
\newcommand{\xequal}[2][]{\ext@arrow 0055{\equalfill@}{#1}{#2}}
\def\equalfill@{\arrowfill@\Relbar\Relbar\Relbar}
\makeatother

\newcommand{\mto}{\mapsto}

\newcommand{\1}{\ensuremath{\ol{\mathrm{P}}}}

\renewcommand{\k}{\ensuremath{\ol{\mathrm{P}}}}

\newcommand{\hou}[3]{{#1}\equiv {#2}\pmod{#3}}

\newcommand{\hl}{\hline}

\newcommand{\h}{\hline}

\renewcommand{\k}[1]{\ensuremath{\left({#1}\right)}}

\newcommand{\re}{\te{Re}\,}

\newcommand{\bca}{\begin{cases}}
\newcommand{\eca}{\end{cases}}

\newcommand{\mug}{\ensuremath{\infty}}

\newcommand{\ff}[2]{\ensuremath{\di\fr{#1}{#2}}}

\newcommand{\bpic}{\begin{picture}}\newcommand{\epic}{\end{picture}}

\newcommand{\beda}{\begin{edaenumerate}}
\newcommand{\eeda}{\end{edaenumerate}}

\newcommand{\cd}{\cdots}

\newcommand{\q}{\quad}

\newcommand{\bq}{\begin{quote}}\newcommand{\eq}{\end{quote}}

\newcommand{\rt}{\sqrt}
\newcommand{\be}{\begin{enumerate}}\newcommand{\ee}{\end{enumerate}}
\newcommand{\bce}{\begin{center}}\newcommand{\ece}{\end{center}}
\newcommand{\bde}{\begin{description}}\newcommand{\ede}{\end{description}}
\newcommand{\bri}{\begin{flushright}}\newcommand{\eri}{\end{flushright}}
\newcommand{\bb}{\begin{block}}\newcommand{\eb}{\end{block}}
\newcommand{\bt}{\begin{thm}}\newcommand{\et}{\end{thm}}
\newcommand{\bpf}{\begin{proof}}\newcommand{\epf}{\end{proof}}
\newcommand{\bex}{\begin{ex}}\newcommand{\eex}{\end{ex}}
\newcommand{\bexr}{\begin{exr}}\newcommand{\eexr}{\end{exr}}
\newcommand{\bft}{\begin{fact}}\newcommand{\eft}{\end{fact}}
\newcommand{\brk}{\begin{rmk}}\newcommand{\erk}{\end{rmk}}
\newcommand{\ba}{\begin{align*}}\newcommand{\ea}{\end{align*}}
\newcommand{\bexe}{\begin{exe}}\newcommand{\eexe}{\end{exe}}
\newcommand{\tn}{\textnormal}

\newcommand{\bit}{\begin{itemize}}\newcommand{\eit}{\end{itemize}}

\newcommand{\bcm}{}
\newcommand{\ol}{\overline}

\newcommand{\fr}{\frac}
\newcommand{\cc}{\ensuremath{\mathbf{C}}}

\newcommand{\bd}{\begin{defn}}\newcommand{\ed}{\end{defn}}
\newcommand{\bp}{\begin{prop}}\newcommand{\ep}{\end{prop}}
\newcommand{\p}{\ensuremath{\pi}}
\newcommand{\eh}{\emph}

\newcommand{\te}{\text}

\newcommand{\di}{\displaystyle}\renewcommand{\a}{\ensuremath{\bm{a}}}

\newcommand{\f}{\frac}

\newcommand{\z}{\ensuremath{\bm{z}}}
\newcommand{\np}{\newpage}

\renewcommand{\a}{\alpha}

\usepackage[dvipdfmx]{graphicx}
\usepackage[dvipsnames]{xcolor}
\usepackage{float,afterpage}
\usepackage{asymptote,layout}
\usepackage{wrapfig,epic}

\usepackage{geometry}
\usepackage{exscale,latexsym,bm}
\usepackage{amssymb,enumerate,amsmath,amsthm,amsfonts}
\usepackage{verbatim,fancybox,wasysym,fancyhdr,type1cm}
\usepackage[frame,all,poly,curve,knot,arrow]{xy}
\usepackage{colortbl}
\usepackage{boxedminipage}
\usepackage{multirow}

\graphicspath{{./figs/}}
\theoremstyle{definition}
\newtheorem{thm}{Theorem}[section]

\newtheorem{prop}[thm]{Proposition}\newtheorem{cor}[thm]{Corollary}

\newtheorem{exr}[thm]{Exercise}

\newtheorem{ex}[thm]{Example}

\newtheorem{defn}[thm]{Definition}\newtheorem{rmk}[thm]{Remark}
\newtheorem{fact}[thm]{Fact}
\newtheorem{block}[thm]{}
\newtheorem*{exe}{Exercise}

\renewcommand{\a}{\alpha}
\renewcommand{\z}{\zeta}
\renewcommand{\tom}{\to \mg}

\begin{document}
\renewcommand{\h}{\hline}
\renewcommand{\arraystretch}{1.5}
\newcommand{\fib}{\text{fib\,}}
\newcommand{\fdp}{\text{FDP\,}}

\title{New relations on zeta and $L$ functions}
\author{Masato Kobayashi}
\date{\today}                                       
\subjclass[2020]{Primary:11M06;\,Secondary:11M41;}
\keywords{Bernoulli numbers, 
Dirichlet $L$ function, 
Euler numbers, 
generating function, Taylor series, 
zeta function.}
\address{Masato Kobayashi\\
Department of Engineering\\
Kanagawa University, 3-27-1 Rokkaku-bashi, Yokohama 221-8686, Japan.}
\email{masato210@gmail.com}


\maketitle
\begin{abstract}
We prove new relations 
on zeta function at even arguments 
and Dirichlet $L$ function at odd. 
The key idea is to make use of Taylor series and partial fraction decomposition of cotangent and secant functions as we discuss in calculus and complex analysis.
\end{abstract}
\tableofcontents

\section{Introduction}
\eh{Riemann zeta function}
\[
\z(s)=
\ff{1}{1^{s}}+\ff{1}{2^{s}}+\ff{1}{3^{s}}+\ff{1}{4^{s}}+\cd, \q \re{(s)}>1
\]
and the following type of 
\eh{Dirichlet $L$ function} (also called 
\eh{Dirichlet beta function}) 
\[
L(s)=
\ff{1}{1^{s}}-\ff{1}{3^{s}}+\ff{1}{5^{s}}-\ff{1}{7^{s}}+\cd
, \q \re{(s)}\ge 1
\]
are of great importance in number theory.

\begin{thm}[Euler]
Set $\z(0)=-\f{1}{2}$. Then for each integer $m\ge 0$,
\begin{align*}
	\z(2m)&=(-1)^{m+1}\ff{(2\p)^{2m}}{2(2m)!}B_{2m},
	\\L(2m+1)&=\k{\ff{\p}{2}}^{2m+1}\ff{1}{2(2m)!}E_{2m}
\end{align*}
where $B_{2m}$, $E_{2m}$ are Bernoulli and Euler numbers as in Tables \ref{bno}, \ref{eno}.
\end{thm}
See Ayoub \cite{ayoub} for history of these series and numbers.

{\renewcommand{\arraystretch}{2}
\begin{table}
\caption{signed Bernoulli numbers 
$(B_{3}=B_{5}=\cd =0$)}
\label{bno}
\begin{center}
\begin{tabular}{c|ccccccccccccccccccccccccc}
$m$&0&1&2&4&6&8\\\h
$B_{m}$&\dep{17pt}1&$- \f12$&$\f{1}{6}$\hgt{20pt}&$-\f{1}{30}$
&$\f{1}{42}$&$-\f{1}{30}$\\\hline
$m$&10&12&14&16&18&$\cd$\\\hline
$B_{m}$&\dep{17pt}$\f{5}{66}$&$-\f{691}{2730}$&$\f{7}{6}$&$-\f{3617}{510}$&
$\f{43867}{798}$
&$\cd$\\\h
\end{tabular}
\end{center}
\end{table}%
}

{\renewcommand{\arraystretch}{1.25}
\begin{table}
\caption{unsigned Euler numbers}
\label{eno}
\begin{center}
\begin{tabular}{c|ccccccccccccccccccc}
$m$&0&1&2&3&4&5&6&7&8&9&$\cd$\\\hline
$E_{m}$&1&0&1&0&5&0&61&0&1385&0&$\cd$\\
\end{tabular}
\end{center}
\end{table}%
}

{\renewcommand{\arraystretch}{2.25}
\begin{table}
\caption{known values of $\z$ and $L$ functions}
\label{tokushu}
\begin{center}
\begin{tabular}{c|ccccccccccccccccccc}
$m$&0&1&2&3&4&5&6&7&8\\\hline
$\zeta(m)$\dep{16pt}&
$-\ff{\,1\,}{2}$&
$+\mug$&$\ff{\p^{2}}{6}$&&$\ff{\p^{4}}{90}$&&$\ff{\p^{6}}{945}$&&
$\ff{\p^{8}}{9450}$
\\\h
$L(m)$\dep{16pt}&&
$\ff{\p}{4}$&&$\ff{\p^{3}}{32}$&
&$\ff{5\p^{5}}{1536}$
&&
$\ff{61\p^{7}}{184320}$
&&\\\h
\end{tabular}
\end{center}
\end{table}%
}

{\renewcommand{\arraystretch}{2.5}
\begin{table}[h!]
\caption{Summary of this article}
\begin{center}
\label{gl}
\begin{tabular}{|c|c|}
\hl $\z$ &$L$  \\\hline
$\z(2m)=(-1)^{m+1}\ff{(2\p)^{2m}}{2(2m)!}B_{2m}$&
$L(2m+1)=\k{\ff{\p}{2}}^{2m+1}\ff{1}{2(2m)!}E_{2m}$
\\\hl
$1<\z(m)<2$ $(m>1)$&$0<L(m)<1$ $(m\ge1)$\\\hl
$\z(m)\to 1$ $(m\tom)$&$L(m)\to 1$ $(m\tom)$\\\hl
$\dsum_{m=2}^{\mug}(\z(m)-1)=1$
&
Theorem \ref{mth1}
\\\hl
$-\ff{\piz}{2}\cot(\piz)=
\dsum_{m=0}^{\mug}\z(2m)z^{2m}$
&
Theorem \ref{mth2}
\\\hl
$\z(2m)=
\ff{2}{2m+1}
\dsum_{k=1}^{m-1}\z(2k)\z(2m-2k)$&
Theorem \ref{mth3}\\\hl
\multicolumn{2}{|c|}{Theorem \ref{mth4} (mixed formula)}
\\\hl
\end{tabular}
\end{center}
\end{table}}

The aim of this article is to find new relations on 
$\{\z(2m)\}_{m=0}^\mg$ and 
$\{L(2m+1)\}_{m=0}^\mg$. Often many researchers discover infinite series involving $\z(2m)$ as follows. 
Borwein-Bradley-Crandall \cite[p.255, 262-263]{borwein}:
\[
\dsum_{k=1}^{m-1}\z(2k)\z(2m-2k)=
\k{m+\f{1}{2}}\z(2m),
\]
\[
\dsum_{m=1}^{\mug} (\z(2m)-1)=
\ff{3}{4},
\]
\[
\dsum_{m=1}^{\mug} \ff{\z(2m)-1}{m}=\log2,
\]
\[
\dsum_{m=1}^{\mug} \ff{\z(2m)-1}{2^{2m}}=\ff{1}{6},
\]
\[
\dsum_{m=1}^{\mug} \ff{\z(2m)-1}{4^{2m}}=\ff{13}{30}
-\ff{\p}{8},
\]
\[
\dsum_{m=1}^{\mug} \ff{\z(2m)-1}{8^{2m}}=\ff{61}{126}
-\ff{\p}{16}\rt{\ff{\rt{2}+1}{\rt{2}-1}},
\]
\[
\dsum_{m=1}^{\mug} (\z(4m)-1)=\ff{7}{8}-\ff{\p}{4}\coth(\p).
\]
Choi \cite[p.388]{choi}:
\[
\dsum_{m=1}^{\mug} \ff{\z(2m)-1}{m+1}=\ff{3}{2}-\log\pi. 
\]
Srivastava \cite[p.133, 136]{srivastava}:
\[
\dsum_{m=1}^{\mug} \ff{\z(2m)}{m(2m+1)}=\log(2\p)-1,
\]
\[
\dsum_{m=1}^{\mug} \ff{\z(2m)}{m\,2^{2m}}=\log\k{\ff{1}{2}\p},
\]
\[
\dsum_{m=1}^{\mug} \ff{\z(2m)}{(2m+1)2^{2m}}=\ff{1}{2}-\log2.
\]
We will see several new series in Example \ref{ex3}.

Presumably, we understand zeta and $L$ functions together. Thus, we will extend known results on $\z(2m)$ to the ones on $L(2m+1)$ as Theorems \ref{mth1}, \ref{mth2} and \ref{mth3}. 
Furthermore, we will prove the equality (which we call the \eh{mixed formula}) involving \eh{both} of $\{\z(2m)\}$ and $\{L(2m+1)\}$ 
as Theorem \ref{mth4}; 
Table \ref{gl} shows summary of this article.



\section{Fractional parts of $\zeta(2m)$ and $L(2m+1)$}
\label{s2}

Throughout 
$j, k, l, m, n$ denote nonnegative integers and $z$ a complex number.

 
Since $\zeta(m)=1+\f{1}{2^{m}}+\cd>1$ and 
\[
2>
\ff{\p^{2}}{6}=\z(2)>\z(3)>\z(4)>\z(5)>\cd,
\]
we have $1<\z(m)<2$ for all $m\ge 2$.
That is, $\z(m)-1$ is the \eh{fractional part} of $\z(m)$. For example, 
\begin{align*}
	\zeta (2)-1&=0.6449\cdots,
	\\\zeta (3)-1&=0.2020\cdots,
	\\\zeta (4)-1&=0.0823\cdots,
	\\\zeta (5)-1&=0.0369\cdots.
\end{align*}
Shallit-Zikan gave the following surprising identity:
\begin{thm}[Shallit-Zikan \cite{shallit}]
\label{sz}
\begin{align*}
	\dsum_{m=2}^{\mug}(\z(m)-1)&=1.
	\end{align*}
\end{thm}
Some other consequences are 
\begin{align*}
	\dsum_{m=1}^{\mug}(\z(2m+1)-1)&=\ff{1}{4},
	\\
	\dsum_{m=1}^{\mug}(\z(2m)-1)&=\ff{3}{4},
	\\
	\dsum_{m=2}^{\mug}(-1)^{m}(\z(m)-1)&=\ff{1}{2}.
\end{align*}
See also Bibiloni-Paradis-Viader \cite{bibiloni} 
and Choi-Quine-Srivastava \cite{choi} 
for Euler-Goldbach Theorem, the origin of such equalities.

Now what can we say about the sum of fractional parts of $\{L(2m+1)\}$? 
As the proposition just below says, 
things are little different.
\begin{prop}
For each $m\ge 1$, we have $0\le L(m)\le 1$.
\end{prop}
\begin{proof}
For fixed $m\ge1$, let 
\[
S_{N}= \sum_{n=1}^{N}\ff{(-1)^{n+1}}{(2n-1)^{m}}
\]
be a partial sum of $L_{}(m)$. It follows that 
\begin{align*}
	S_{2N+1}&=\sum_{n=1}^{2N+1}\ff{(-1)^{n+1}}{(2n-1)^{m}}
	\\&=1-\ff{1}{3^{m}}+\ff{1}{5^{m}}-\cd -
	\ff{1}{(4N-1)^{m}}+\ff{1}{(4N+1)^{m}}
	\\&=1+
	\underbrace{\k{-\ff{1}{3^{m}}+\ff{1}{5^{m}}}}_{<0}
+\cd 
+
\underbrace{\k{-\ff{1}{(4N-1)^{m}}+\ff{1}{(4N+1)^{m}}}}_{<0}
	\\&<1
\end{align*}
and 
\[
S_{2N+2}=S_{2N+1}-\ff{1}{(4N+3)^{m}}<S_{2N+1}<1.
\]
Hence $L_{}(m)= \di\lim_{N\to \mug}S_{N}\le 1$.
Moreover, there exists the following infinite product 
(known as \eh{Euler product})
\[
L(m)=\di \prod_{\hou{p}{1}{4}} (1-p^{-m})^{-1}
\di \prod_{\hou{p}{3}{4}}
 (1+p^{-m})^{-1}
\]
with $p$ prime numbers and each factor $(1\pm p^{-m})^{-1}\ge 0$. Thus, $L(m)\ge 0$ and this completes the proof.
\end{proof}
For this reason, let us consider $1-L(m)$ instead of 
the fractional part of $L(m)$. For example, 
\begin{align*}
	1-L(1)&=0.21460\cd,
	\\1-L(2)&=0.08403\cd,
	\\1-L(3)&=0.03105\cd,
	\\1-L(4)&=0.01105\cd.
\end{align*}
The following is an analogy of 
Theorem \ref{sz}.
\begin{thm}\label{mth1}
\begin{eqnarray}
	\dsum_{m=1}^{\mug}(1-L(m))=\ff{1}{2}\log2=
0.34657\cd,
\label{eq1}
\end{eqnarray}
\begin{eqnarray}
\dsum_{m=1}^{\mug}(1-L(2m))=\ff{1}{2}\log2-\ff{1}{4},
\label{eq2}
\end{eqnarray}
\begin{eqnarray}
\dsum_{m=1}^{\mug}(1-L(2m-1))=\ff{1}{4},
\label{eq3}
\end{eqnarray}
\begin{eqnarray}
\dsum_{m=1}^{\mug}(-1)^{m-1}(1-L(m))=\ff{1}{2}(\log2-1).
\label{eq4}
\end{eqnarray}

\end{thm}
\begin{proof}
First, we remark that 
\[
\dsum_{m=2}^{\mug} 
\dsum_{n=2}^{\mug} \abs{\ff{(-1)^{n}}{(2n-1)^{m}}}
=
\dsum_{m=2}^{\mug} 
\dsum_{n=2}^{\mug} \ff{1}{(2n-1)^{m}}
\le \dsum_{m=2}^{\mug} (\z(m)-1)
=1
\]
(absolutely convegent)
and then 
\[
\dsum_{m=1}^{\mug} (1-L(m))
=(1-L(1))
+
\dsum_{m=2}^{\mug} 
\dsum_{n=2}^{\mug} 
\ff{(-1)^{n}}{(2n-1)^{m}}
\]
must be convergent;
in particular, we can exchange summations
(likewise below).
	
\begin{align*}
	\dsum_{m=1}^{\mug} (1-L(m))&=
	\dsum_{m=1}^{\mug} 
	\dsum_{n=2}^{\mug} \ff{(-1)^{n}}{(2n-1)^{m}}
	\\&=
	\dsum_{n=2}^{\mug} \ff{(-1)^{n}}{2n-1}
	\dsum_{m=0}^{\mug}\k{\ff{1}{2n-1}}^{m}
	\\&=
	\dsum_{n=2}^{\mug} \ff{(-1)^{n}}{2n-1}
	\k{\ff{1}{1-\f{1}{2n-1}}}
	\\&=
	\ff{1}{2}
	\dsum_{n=2}^{\mug} \ff{(-1)^{n}}{n-1}=
	\ff{1}{2}\log2.
\end{align*}

\begin{align*}
	\dsum_{m=1}^{\mug} (1-L(2m))&=
	\dsum_{m=1}^{\mug} \dsum_{n=2}^{\mug} \ff{(-1)^{n}}{(2n-1)^{2m}}
	\\&=
	\dsum_{n=2}^{\mug} 
	\ff{(-1)^{n}}{(2n-1)^{2}}
	\dsum_{m=0}^{\mug}\k{\ff{1}{2n-1}}^{2m}
	\\&=
	\dsum_{n=2}^{\mug} 
	\ff{(-1)^{n}}{(2n-1)^{2}}
	\k{\ff{1}{1-
	\k{\f{1}{2n-1}}^{2}
	}}
	\\&=
	\dsum_{n=2}^{\mug} \ff{(-1)^{n}}{(2n-2)2n}
	\\&=
	\ff{1}{4}
	\k{\dsum_{n=2}^{\mug} \ff{(-1)^{n}}{n-1}
	-
	\dsum_{n=2}^{\mug} \ff{(-1)^{n}}{n}}
	\\&=\ff{1}{4}(\log2+(\log2-1))
	\\&=\ff{1}{2}\log2-\ff{1}{4}.
\end{align*}
Moreover, (\ref{eq3})
 is $(\ref{eq1})-(\ref{eq2})$ and 
 (\ref{eq4}) is 
 $(\ref{eq2})-
 (\ref{eq3})$.
\end{proof}


\section{Generating function}
\label{s3}

To study sequences $\{\z(2m)\}$ and $\{L(2m+1)\}$, 
it is helpful to find their generating functions. In the sequel, we will often use the following fact in complex analysis implicitly.
\begin{fact}
If $F(z)=
\sum_{m=0}^{\mug}\a_{m}z^{m}$ ($\a_{m}\in\cc$) is a convergent power series with the radius of convergence $R$, then so is $F'(z)$ and moreover it is given by the power series 
$\sum_{m=0}^{\mug}m\a_{m}z^{m-1}.$
\end{fact}


\begin{prop}\label{p1}
For $|z|<1$, 
\[
-\ff{\pi z}{2}\cot(\pi z)=
\dsum_{m=0}^{\mug}
\z(2m)z^{2m}.
\]
\end{prop}

\begin{proof}Recall from complex analysis that 
\[
\cot(\pz)=
\ff{1}{\pz}+\ff{1}{\p} \dsum_{n=1}^{\mug}\ff{2z}{z^{2}-n^{2}}.
\]
Now let $|z|<1$. Then, we have $\abs{\f{z}{n}}^{2}<1$ for all $n\ge 1$ so that 
\begin{align*}
	-\ff{\pz}{2}\cot(\pz)&=-\ff{1}{2}-
	\dsum_{n=1}^{\mug}\ff{z^{2}}{z^{2}-n^{2}}
	\\&=-\ff{1}{2}+
	\dsum_{n=1}^{\mug}\ff{z^{2}}{n^{2}\k{1-\k{\f{z}{n}}^{2}
}
}
	\\&=-\ff{1}{2}+ \dsum_{n=1}^{\mug}\ff{1}{n^{2}}
	\k{\dsum_{m=0}^{\mug}\k{\ff{z}{n}}^{2m}}z^{2}
	\\&=-\ff{1}{2}+
	\dsum_{m=0}^{\mug}
	\k{\dsum_{n=1}^{\mug}
	\ff{1}{n^{2m+2}}}
	z^{2m+2}
	\\&=\z(0)+\dsum_{m=0}^{\mug}\z(2m+2)z^{2m+2}
	\\&=\dsum_{m=0}^{\mug}
\z(2m)z^{2m}.
\end{align*}
\end{proof}
\begin{ex}\label{ex3}
Set $f(z)=-\f{\pi z}{2}\cot(\pi z)$ as above.
Observe that 
\begin{align*}
	\dsum_{m=0}^{\mug} \ff{\z(2m)}{2^{2m}}&=0,
	\\\dsum_{m=0}^{\mug} \ff{\z(2m)}{3^{2m}}&=-\ff{\p}{6\rt{3}},
	\\\dsum_{m=0}^{\mug} \ff{\z(2m)}{4^{2m}}&=-\ff{\p}{8},
	\\\dsum_{m=0}^{\mug} \ff{\z(2m)}{5^{2m}}&=
-\ff{\rt{2}}{10\rt{5-\rt{5}}}\k{\ff{1+\rt{5}}{2}\p},
	\\\dsum_{m=0}^{\mug} \ff{\z(2m)}{6^{2m}}&=-\ff{\rt{3}}{12}\p
\end{align*}
\end{ex}
are $f(\f{1}{2})$, $f(\f{1}{3})$, $f(\f{1}{4})$, $f(\f{1}{5})$, $f(\f{1}{6})$, respectively 
$\k{\text{recall }\cot\f{\p}{5}=
\f{\f{1+\rt{5}}{4}}{\rt{\f{5-\rt{5}}{8}}}}$.
There are some consequences. We see that 
\[
zf'(z)
=\dsum_{m=0}^{\mug}2m\z(2m)z^{2m},
\]
that is, 
\[
\dsum_{m=0}^{\mug}2m\z(2m)z^{2m}=
-\ff{\piz}{2}\k{
\cot\piz+\piz(-\cot^{2}(\piz)-1
)}
.
\]
Then observe also that 
\begin{align*}
\dsum_{m=0}^{\mug}\ff{2m}{2^{2m}}\z(2m)	
&
=
\ff{1}{2}f'\k{\f{1}{2}}
=
	\ff{\p^{2}}{8},
	\\
	\dsum_{m=0}^{\mug}\ff{2m}{3^{2m}}\z(2m)&=\f{1}{3}f'\k{\f{1}{3}}
	=
	-\ff{\p}{6}\k{\ff{1}{\rt{3}}-\ff{4}{9}\p}
\end{align*}
and so on.

\begin{thm}\label{mth2}
For $|z|<1$,
\[
\ff{\pi z}{4}\sec\ff{\pi z}{2}=
\dsum_{m=0}^{\mug}L(2m+1)z^{2m+1}.
\]
\end{thm}
\begin{proof}
Recall from complex analysis that 
\[
\sec z=\dsum_{n=1}^{\mug}\ff{(-1)^{n}(2n-1)\p}{z^{2}-\k{\f{2n-1}{2}\p}^{2}}.\]
Now let $|z|<1$. 
Then, for all $n\ge 1$, we have $\abs{\f{z}{2n-1}}^{2}<1$ so that 
\begin{align*}
	\ff{\pi z}{4}\sec\ff{\pi z}{2}&=
	\ff{\pi z}{4}
\dsum_{n=1}^{\mug}\ff{(-1)^{n}(2n-1)\p}{
\k{\f{\piz}{2}}^{2}-
\k{\f{2n-1}{2}\p}^{2}
}
	\\&=z
	\dsum_{n=1}^{\mug}\ff{(-1)^{n+1}(2n-1)}{
	(2n-1)^{2}
	\k{1-\k{\f{z}{2n-1}}^{2}
}
	}
	\\&=z
	\dsum_{n=1}^{\mug}\ff{(-1)^{n+1}}{2n-1}
	\dsum_{m=0}^{\mug}\k{\ff{z}{2n-1}}^{2m}
	\\&=
	\dsum_{m=0}^{\mug}
	\k{\dsum_{n=1}^{\mug}\ff{(-1)^{n+1}}{(2n-1)^{2m+1}}}z^{2m+1}
	\\&=\dsum_{m=0}^{\mug}L(2m+1)z^{2m+1}.
\end{align*}
\end{proof}

\begin{cor}
\begin{align*}
	\dsum_{m=0}^{\mug}\ff{L(2m+1)}{2^{2m+1}}&=
	\ff{\rt{2}}{8}\p,
	\\\dsum_{m=0}^{\mug}\ff{L(2m+1)}{3^{2m+1}}&=\ff{\p}{6\rt{3}},
	\\\dsum_{m=0}^{\mug}\ff{L(2m+1)}{4^{2m+1}}&=\ff{\p}{8\rt{2+\rt{2}}},
	\\
	\dsum_{m=0}^{\mug}\ff{L(2m+1)}{5^{2m+1}}&=
	\ff{\rt{2}}{10\rt{5+\rt{5}}}\p,
	\\\dsum_{m=0}^{\mug}\ff{L(2m+1)}{6^{2m+1}}&=
	\ff{\rt{6}-\rt{2}}{24}\p.
\end{align*}
%
\end{cor}
\begin{proof}
Set $g(z)=\f{\pi z}{4}\sec\f{\pi z}{2}$ as 
above. These identities are $g\k{\f{1}{2}}$, $g\k{\f{1}{3}}$, $g\k{\f{1}{4}}$, $g\k{\f{1}{5}}$, $g\k{\f{1}{6}}$ with 
\[
\sec\ff{\p}{8}=\ff{2}{\rt{2+\rt{2}}}, \q
\sec\ff{\p}{10}=\ff{2\rt{2}}{\rt{5+\rt{5}}},\q
\sec\ff{\p}{12}=\rt{6}-\rt{2}.
\]
Moreover, 
\[
zg'(z)=\ff{\piz}{4}
\sec\ff{\piz}{2}\k{1+\ff{\pi z}{2}\tan\ff{\piz}{2}}
\]
leads us to 
\begin{align*}
	\dsum_{m=0}^{\mug}\ff{2m+1}{2^{2m+1}}
L(2m+1)	
&=\ff{1}{2}g'\k{\f{\,1\,}{2}}
	=
\ff{\rt{2}}{8}\p\k{1+\ff{\p}{4}},
	\\
	\dsum_{m=0}^{\mug}\ff{2m+1}{3^{2m+1}}L(2m+1)
	&=\ff{1}{3}g'\k{\f{\,1\,}{3}}
	=
\ff{1}{6\rt{3}}\p\k{1+\ff{\p}{18\rt{3}}}
\end{align*}
and so on.
\end{proof}

\section{Convolution}
\label{s4}

It is well-known that $\{\z(2m)\}$ satisfies the following relation.
\begin{fact}[{\cite[p.255]{borwein}}]
For $m\ge2$,
\[
\z(2m)=\ff{2}{2m+1}
\sum_{
\substack{j, k>0\\j+k=m
}
}
\z(2j)\z(2k).
\]
\end{fact}
The following is an analogous result of this.

\begin{thm}\label{mth3}
For $m\ge 1$, we have 
\[
L(2m+1)=
\ff{1}{(2m-1)2m}
\k{8
\sum_{
\substack{j, k, l\ge 0\\j+k+l=m-1
}
}L(2j+1)L(2k+1)L(2l+1)
-
\f{\p^{2}}{4}L(2m-1)}.
\]
\end{thm}

\begin{proof}
Note that
\begin{eqnarray}
\k{\sec\ff{\piz}{2}}''=
\ff{\p^{2}}{4}
\k{2\sec^{3}\ff{\piz}{2}-\sec\ff{\piz}{2}
}.
\label{eq5}
\end{eqnarray}
Let us now express each term in both sides as power series:
\[
\k{\sec\ff{\piz}{2}}''=
\ff{4}{\p}
\dsum_{m=0}^{\mug}(2m)(2m-1)L(2m+1)z^{2m-2},
\]
\[
\ff{\p^{2}}{4}
\k{2\sec^{3}\ff{\piz}{2}}
=
\ff{\p^{2}}{2}
\k{
\ff{4}{\p}
\dsum_{m=0}^{\mug}L(2m+1)z^{2m}
}^{3}
\]
\[
=
\ff{32}{\p}
\dsum_{m=0}^{\mug}
\k{\sum_{
\substack{j, k, l\ge 0\\j+k+l=m
}
}L(2j+1)L(2k+1)L(2l+1)}z^{2m},
\]
\[
-\ff{\p^{2}}{4}\sec\ff{\piz}{2}
=-\pi \dsum_{m=0}^{\mug}L(2m+1)z^{2m}.
\]
Equate the coefficients of $z^{2(m-1)}$ $(m\ge1)$ 
in (\ref{eq5}):
\[
\ff{4}{\p}(2m)(2m-1)L(2m+1)
=
\ff{32}{\p}
\k{\sum_{
\substack{j, k, l\ge 0\\j+k+l=m-1
}
}L(2j+1)L(2k+1)L(2l+1)}
-\p L(2m-1).
\]
Conclude that 
\[
L(2m+1)=
\ff{1}{(2m-1)2m}
\k{8
\sum_{
\substack{j, k, l\ge 0\\j+k+l=m-1
}
}L(2j+1)L(2k+1)L(2l+1)
-
\f{\p^{2}}{4}L(2m-1)}.
\]
\end{proof}
For example, 
\[
L(5)=\ff{1}{3\cdot 4}\k{8\cdot 3L(1)^{2}L(3)-
\ff{\p^{2}}{4}L(3)}
=\ff{1}{12}
\k{24
\k{\ff{\p}{4}}^{2}
-\ff{\p^{2}}{4}}
\ff{1}{32}\p^{3}
=
\ff{5}{1536}\p^{5}.
\]

\section{Mixed formula}
\label{s5}
We have seen many results for each of 
zeta function and $L$ function. 
Here, we prove the identity involving 
\eh{both} of $\{\z(2m)\}$ and $\{L(2m+1)\}$.
\begin{thm}[mixed formula]\label{mth4}
For $m\ge1$,
\[
L(2m+1)=
\ff{1}{4m}
\k{
\sum_{
{\substack{j, k>0\\j+k=m
}}
}
\ff{\z(2j)}{2^{2j}}8kL(2k+1)+\ff{\p^{2}}{2}L(2m-1)
}.
\]
\end{thm}
\begin{proof}
Note that 
\[
\k{\sec\ff{\piz}{2}}'=
\ff{\p}{2}\tan\ff{\piz}{2}\sec\ff{\piz}{2}
\]
and so 
\[
\k{-\ff{\piz}{4}\cot\ff{\piz}{2}}
\k{\sec\ff{\piz}{2}}'
=
-\ff{\pi^{2}z}{8}\sec\ff{\piz}{2},
\]
that is, 
\[
\k{\dsum_{j=0}^{\mug}\z(2j)\k{\ff{z}{2}}^{2j}}
\k{\ff{4}{\p}
\dsum_{k=0}^{\mug}2kL(2k+1)z^{2k-1}}
=
-\ff{\p}{2}
\k{
\dsum_{m=0}^{\mug}L(2m+1)z^{2m+1}
}.
\]
Now equate the coefficients of $z^{2m+1}$ $(m\ge 0)$ in both sides to obtain 
\[
\sum_{
\substack{j, k\ge 0\\j+k=m+1
}
}\ff{\z(2j)}{2^{2j}}\ff{4}{\p}2kL(2k+1)=-\ff{\p}{2}L(2m+1).
\]
In the sum on the left hand side, 
the term for $(j, k)=(0, m+1)$ is 
\[
\ff{\z(0)}{2^{0}}\ff{\p}{4}2(m+1)L(2(m+1)+1)=
-\ff{4(m+1)}{\p}L(2m+3)
\]
and for $(j, k)=(m+1, 0)$ is 0. 
Therefore, we have 
\[
-\ff{4(m+1)}{\p}L(2m+3)+
\sum_{
\substack{j, k>0\\j+k=m+1
}
}\ff{\z(2j)}{2^{2j}}\ff{4}{\p}2kL(2k+1)
=
-\ff{\p}{2}L(2m+1)
\]
and so 
\[
L(2m+3)=
\ff{\p}{4(m+1)}
\k{
\sum_{
\substack{j, k>0\\j+k=m+1
}
}\ff{\z(2j)}{2^{2j}}\ff{8k}{\p}L(2k+1)
+
\ff{\p}{2}L(2m+1)}.
\]
With $m\mto m-1$, we conclude that 
\[
L(2m+1)=
\ff{1}{4m}
\k{
\sum_{
\substack{j, k>0\\j+k=m
}
}
\ff{\z(2j)}{2^{2j}}8kL(2k+1)+\ff{\p^{2}}{2}L(2m-1)
}.
\]
\end{proof}
For example, let $m=3$.
\[
L(7)=
\ff{1}{12}
\k{
\ff{\z(2)}{2^{2}}16L(5)+\ff{\z(4)}{2^{4}}8L(3)
+\ff{\p^{2}}{2}L(5)
}
\]
\[
=
\ff{1}{12}
\k{
4
\ff{\p^{2}}{6}
\ff{5\p^{5}}{1536}
+
\f{1}{2}
\ff{\p^{4}}{90}
\ff{\p^{3}}{32}
+
\ff{\p^{2}}{2}
\ff{5\p^{5}}{1536}
}=
\ff{61}{184320}\p^{7}.\]

\section{Conclusion}
In this note, we proved many identities on $\{\z(2m)\}_{m=0}^{\mg}$
 and $\{L(2m+1)\}_{m=0}^{\mg}$. The key idea is to 
 make use of Taylor series and partial fraction decomposition of cotangent and secant functions as we sometimes discuss in calculus and complex analysis.
 Indeed, this method is powerful enough to 
 find more formulas.
For example, we can derive other relations from 
\[
2\cot(z)\cot(2z)=\cot^{2}z-1, 
\]
\[
1+\cot^{2}z=\cot^{2}z\sec^{2}z,
\]
\[
-\ff{\piz}{4}
(\cot(\piz)+\coth(\piz))=
\dsum_{m=0}^{\mug}\z(4m)z^{4m}
,
\]
\[
-\ff{\piz}{4}
(\cot(\piz)-\coth(\piz))
=
\dsum_{m=0}^{\mug}\z(4m+2)z^{4m+2},
\]
\[
\ff{\piz}{8}\k{\sec\ff{\piz}{2}+\sech\ff{\piz}{2}}
	=
\dsum_{m=0}^{\mug}L(4m+1)z^{4m+1},
\]
\[
\ff{\piz}{8}\k{\sec\ff{\piz}{2}-\sech\ff{\piz}{2}}
=
\dsum_{m=0}^{\mug}L(4m+3)z^{4m+3}
.
\]

It is also easy to translate our results into relations 
for Bernoulli and Euler numbers as Dilcher discussed \cite{dilcher}.

\begin{center}
Acknowledgment.
\end{center}
\begin{quote}
This research arose from Iitaka online seminar in 2020-2021. 
The author would like to thank the organizer 
Shigeru Iitaka and Kouichi Nakagawa for fruitful discussion. He also thanks Satomi Abe and 
Michihito Tobe for supporting his writing the manuscript.
\end{quote}

\end{document}